\documentclass[a4paper,10pt]{elsart}

\usepackage{amsfonts,amssymb,amsmath,xspace}
\usepackage{graphics,pst-tree}
\usepackage{algorithmic,color}
\usepackage{algorithm,natbib}

\newcommand{\C}{\ensuremath{\mathcal C}\xspace}

\begin{document}

\begin{frontmatter}
\title{A variant of the tandem duplication - random loss model of genome rearrangement}
\author[liafa]{Mathilde Bouvel}
\author[liafa]{Dominique Rossin}
\address[liafa]{LIAFA, Universit\'e Paris Diderot, CNRS, Case 7014, 75205 Paris Cedex 13}

\begin{abstract}
In \cite{CCMR06}, Chaudhuri, Chen, Mihaescu and Rao study algorithmic properties of the \emph{tandem duplication - random loss model} of genome rearrangement, well-known in evolutionary biology. In their model, the cost of one step of duplication-loss of width $k$ is $\alpha^k$ for $\alpha =1$ or $\alpha \geq 2$. In this paper, we study a variant of this model, where the cost of one step of width $k$ is 1 if $k \leq K$ and $\infty$ if $k > K$, for any value of the parameter $K \in \mathbb{N} \cup \{\infty\}$. We first show that permutations obtained after $p$ steps of width $K$ define classes of pattern-avoiding permutations. We also compute the numbers of duplication-loss steps of width $K$ necessary and sufficient to obtain any permutation of $S_n$, in the worst case and on average. In this second part, we may also consider the case $K=K(n)$, a function of the size $n$ of the permutation on which the duplication-loss operations are performed.
\end{abstract}

\begin{keyword}
Sorting \sep Permutations \sep Pattern
\PACS 
\end{keyword}

\end{frontmatter}

\section{Introduction}

\subsection{The model}

In the usual models of genome rearrangement, duplications and losses of genes are not taken into account. There were attempts to incorporate them to the classical models, but the consecutive combinatorial complexity of the models so obtained made their study quite difficult. Following \cite{CCMR06}, we focus on the duplication-loss problem by considering the \emph{tandem duplication - random loss model} of genome rearrangement in which genomes are modified \emph{only} by duplications and losses of genes.

One \emph{step} of tandem duplication - random loss, or duplication-loss for short, consists in $(1)$ the tandem duplication of a contiguous fragment of the genome, \emph{i.e.,} the duplicated fragment is inserted immediately after the original fragment, and $(2)$ the loss of one of the two copies of every duplicated gene. We assume that the loss occurs immediately after the duplication of genes, which is, on an evolutionary time-scale, a good approximation to reality. The \emph{width} of a step is the number of duplicated genes. See Figure \ref{ex:duplication-loss} for an example.

\begin{figure}[ht]
\begin{center}
\begin{eqnarray*}
1 \ 2 \ \overbrace{3 \ 4 \ 5 \ 6} \ 7\ & \rightsquigarrow & 1\ 2\ \overbrace{3\ 4\ 5\ 6}\ \overbrace{3 \ 4\ 5\ 6}\ 7  \\
& & \textrm{ (tandem duplication) } \\
 & \rightsquigarrow & 1\ 2\ 3 \ \hspace{-1em} \diagup \, 4\ 5\ 6\ \hspace{-1em} \diagup \, 3 \ 4\ \hspace{-1em} \diagup \, 5\ \hspace{-1em} \diagup \, 6\ 7 \\
& &  \textrm{ (random loss) } \\
& \rightsquigarrow & 1\ 2\ 4\ 5\ 3\ 6\ 7
\end{eqnarray*}
\caption{Example of one step of tandem duplication - random loss of width $4$ \label{ex:duplication-loss}}
\end{center}
\end{figure}

From a formal point of view, a genome consisting of $n$ genes is modelled by a permutation $\pi \in S_n$ of the set of integers $\{1, 2, \ldots, n\}$. In \cite{CCMR06}, the authors define the cost of a duplication-loss step of width $k$ to be $\alpha^k$, $\alpha \geq 1$ being a constant parameter. They suggest that other cost functions can be considered, and in particular affine functions. In this paper, we consider a \emph{piecewise constant} cost function: the cost of a step of width $k$ is $1$ if $k \leq K$ and is infinite for $k > K$, for some fixed parameter $K \in \mathbb{N} \cup \{\infty\}$. Obviously, for this model to be meaningful, we assume that $K \geq 2$. We also consider the possibility that $K = K(n)$ is dependent on the size $n$ of the permutation on which the duplication-loss operations are performed. Both models are generalizations of the \emph{whole genome duplication - random loss model}: it corresponds to the case $\alpha = 1$ in the model of \cite{CCMR06}, $K = \infty$ or $K = K(n) =n$ in our model.

Many models of evolution of permutations are inspired by computational biology
issues: see \cite{BBCP07}, \cite{CL04}, \cite{conf/wabi/Labarre05}, \cite{journals/tcbb/Labarre06} for examples in the literature.

Our model of evolution of permutations can be viewed in the framework of \emph{permuting machines} defined in \cite{AAA+04}. Such a machine takes a permutation in input, and transforms it into an output permutation, the transformation being subject to satisfy the two properties of 
independence with respect to the values and of stability with respect to pattern-involvement (see \cite{AAA+04} for more details). The important point is that the duplication-loss transformation satisfies these two properties. Thus, one duplication-loss step (in one of the models defined above) corresponds to running an adequate permuting machine once. When we will consider permutations obtained after a sequence of duplication-loss steps, it will correspond to permutations obtained in the output of a combination in series of identical permuting machines.

\begin{figure}[h!]
\begin{center}
\psset{unit=0.3cm}
\begin{pspicture}(-4,0)(11,11)
\psgrid[subgriddiv=1,gridwidth=.2pt,griddots=5,gridlabels=0pt](0,0)(8,8)
\rput(0.5,-0.5){{\tiny $1$}}
\rput(1.5,-0.5){{\tiny $2$}}
\rput(2.5,-0.5){{\tiny $3$}}
\rput(3.5,-0.5){{\tiny $4$}}
\rput(4.5,-0.5){{\tiny $5$}}
\rput(5.5,-0.5){{\tiny $6$}}
\rput(6.5,-0.5){{\tiny $7$}}
\rput(7.5,-0.5){{\tiny $8$}}
\rput(-0.5,0.5){{\tiny $1$}}
\rput(-0.5,1.5){{\tiny $2$}}
\rput(-0.5,2.5){{\tiny $3$}}
\rput(-0.5,3.5){{\tiny $4$}}
\rput(-0.5,4.5){{\tiny $5$}}
\rput(-0.5,5.5){{\tiny $6$}}
\rput(-0.5,6.5){{\tiny $7$}}
\rput(-0.5,7.5){{\tiny $8$}}
\pscircle*(0.5,5.5){0.3}
\pscircle*(1.5,7.5){0.3}
\pscircle*(2.5,0.5){0.3}
\pscircle*(3.5,2.5){0.3}
\pscircle*(4.5,4.5){0.3}
\pscircle*(5.5,3.5){0.3}
\pscircle*(6.5,1.5){0.3}
\pscircle*(7.5,6.5){0.3}
\end{pspicture}
\caption{The graphical representation of $\sigma = 68135427$ \label{fig:graphical-representation}}
\end{center}
\end{figure}
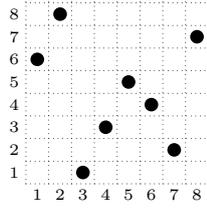

For ease of exposition in some proofs, we will sometimes use a graphical representation of permutations, as shown in Figure \ref{fig:graphical-representation}.

\subsection{Pattern-avoiding classes of permutations}

Though not appearing clearly for the moment, there exist strong links between the duplication-loss model and some pattern-avoiding classes of permutations. Hence, we need to recall a few definitions concerning those classes.

A permutation $\sigma \in S_n$  is a bijective map from $[1..n]$ to itself. The integer $n$ is called the \emph{size} of $\sigma$, denoted $|\sigma|$. We denote by $\sigma_i$ the image of $i$ under $\sigma$. A permutation can be seen as a word $\sigma_1 \sigma_2 \ldots \sigma_n$ containing exactly once each letter $i \in [1..n]$. For each entry $\sigma_i$ of a permutation $\sigma$, we call $i$ its \emph{position} and $\sigma_i$ its \emph{value}.

\begin{defn}
A permutation $\pi \in S_k$ is a \emph{pattern} of a permutation $\sigma \in S_n$ if there is a subsequence of $\sigma$ which is order-isomorphic to $\pi$; in other words, if there is a subsequence $\sigma_{i_1} \sigma_{i_2} \ldots \sigma_{i_k}$ of $\sigma$ (with $1 \leq i_1 < i_2 <\ldots<i_k \leq n$) such that $\sigma_{i_{\ell}} < \sigma_{i_m}$ whenever $\pi_{\ell} < \pi_{m}$. \\
We also say that $\pi$ is \emph{involved} in $\sigma$ and call $\sigma_{i_1} \sigma_{i_2} \ldots \sigma_{i_k}$ an \emph{occurrence} of $\pi$ in $\sigma$. \label{def:pattern}
\end{defn}

We write $\pi \prec \sigma$ to denote that $\pi$ is a pattern of $\sigma$.

A permutation $\sigma$ that does not contain $\pi$ as a pattern is said to \emph{avoid} $\pi$. The class of all permutations avoiding the patterns $\pi_1, \pi_2 \ldots \pi_k$ is denoted $S(\pi_1, \pi_2, \ldots, \pi_k)$, and $S_n(\pi_1, \pi_2, \ldots, \pi_k)$ denotes the set of permutations of size $n$ avoiding $\pi_1, \pi_2, \ldots, \pi_k$. We say that $S(\pi_1, \pi_2, \ldots, \pi_k)$ is a class of pattern-avoiding permutations of \emph{basis} $\{\pi_1, \pi_2, \ldots, \pi_k\}$.

\begin{exmp}
For example $\sigma=1  4  2  5  6  3$ contains the pattern 
$1  3  4  2$, and $1  5  6  3$, $1  4  6  3$, $2 5 6 3$ and 
$1 4 5 3$ are the occurrences of this pattern in $\sigma$. 
But $\sigma \in S(3  2  1)$: $\sigma$ avoids the pattern 
$3  2  1$ as no subsequence of size $3$ of $\sigma$  is isomorphic 
to $3  2  1$, \textit{i.e.}, is decreasing. \label{ex:pattern}
\end{exmp}

\subsection{Outline of the paper}

In the tandem duplication - random loss model described above, we will focus on two kinds of problems. First, as hinted before, we will consider permutations obtained after a certain number of duplication-loss steps, that is to say permutations in output of a combination in series of a certain number of permuting machines. For this, we define the class $\C(K,p)$ as follows: 
\begin{defn}
The class $\C(K,p)$ denotes the class of all permutations obtained from $12\ldots n$ (for any $n$) after $p$ duplication-loss steps of width at most $K$, for some constant parameters $p$ and $K$.
\end{defn}

We do not consider the case $K=K(n)$ here.

Be careful that the duplication-loss steps are not reversible, as noticed in \cite{CCMR06}, and that consequently $\C(K,p)$ is \emph{not} the class of permutations that can be \emph{sorted to} $12\ldots n$ in $p$ steps of duplication-loss of width at most $K$.

Like for the various classes of permutations obtained after a combination in
series of permuting machines considered in \cite{AAA+04}, we obtained
combinatorial properties of $\C(K,p)$ in terms of pattern-avoidance. Namely, we
show that $\C(K,p)$ is a class of pattern-avoiding permutations. In the case
$p=1$ (Section \ref{section:one-step}), we give a precise description of the
basis $B$ of excluded patterns: $B = \{321, 3142, 2143\} \cup D$, $D$ being the
set of all permutations of $S_{K+1}$ that do not start with $1$ nor end with
$K+1$, and containing exactly one descent. In particular, $B$ is of cardinality
$3+2^{K-1}$ and contains patterns of size at most $K+1$. For the general case
(Section \ref{section:p-steps}), we cannot get such a precise result but only a
bound on the size of the excluded patterns: we show that $\C(K,p)$ is a class of
pattern-avoiding permutations whose basis contains patterns of size at most
$(Kp+2)^2-2$.

A second point of view is to examine how many steps of a given width are necessary to obtain any permutation of $S_n$ starting from $12\ldots n$. Namely in Section \ref{section:algo} we fix a width $K$ (constant, or $K=K(n)$) and a size $n$ and search for the number $p$ such that any permutation of $S_n$ can be obtained from $12\ldots n$ in at most $p$ duplication-loss steps of width at most $K$. We describe an algorithm computing a possible scenario of duplications and losses for any $\pi \in S_n$, this scenario involving $\Theta(\frac{n}{K} \log K + \frac{n^2}{K^2})$ duplication-loss steps in the worst case and on average. We also show that $\Omega(\log n + \frac{n^2}{K^2})$ steps are necessary (in the worst case and on average) to obtain any permutation of $S_n$ from $12\ldots n$. These upper and lower bounds coincide in most cases.

\section{Characterization with excluded patterns}

Before focusing on the classes $\C(K,1)$ and $\C(K,p)$ defined for our model, we will get back to the simpler whole genome duplication - random loss model (corresponding to $K=\infty$ in our model, but defined previously by other authors). We will not prove new theorems, but will interprete the existing results from the pattern-avoidance point of view.

\subsection{The whole genome duplication - random loss model through the pattern-avoidance prism}
\label{section:whole-genome}

Let us recall that in the whole genome duplication - random loss model, any duplication-loss step has cost $1$, so that we can consider \emph{w.l.o.g} that the duplicated fragment is the whole permutation at any step. The cost of obtaining a permutation $\sigma \in S_n$ from the identity is just the minimal number of steps of a duplication-loss scenario transforming $12\ldots n$ into $\sigma$.

A statistics of permutations that matters for our purpose is their number of \emph{descents}.

\begin{defn}
Given a permutation $\sigma$ of size $n$, we say that there is a \emph{descent} (resp. \emph{ascent}) at position $i$, $1 \leq i \leq n-1$, if $\sigma_i > \sigma_{i+1}$ (resp. $\sigma_i < \sigma_{i+1}$ ). We write $desc(\sigma)$ the number of descents of the permutation $\sigma$. \label{def:descent}
\end{defn}
\begin{exmp}
For example, $\sigma = 5 2 4 3 1 6$ has $3$ descents, namely at positions $1$, $3$ and $4$. \label{ex:descent}
\end{exmp}

A permutation $\sigma$ of size $n$ has at most $n-1$ descents, the case of $n-1$ descents exactly corresponding to the reversed identity permutation $n(n-1) \ldots 21$. It is also of common knowledge that the average number of descents among permutations of size $n$ is $\frac{n-1}{2}$.

In \cite{CCMR06}, the authors prove the following theorem.
\begin{thm}
Let $\sigma \in S_n$. In the whole genome duplication - random loss model, $\lceil \log_2 (desc(\sigma) + 1) \rceil$ steps are necessary and sufficient to obtain $\sigma$ from $12\ldots n$.
\label{thm:whole-genome}
\end{thm}

It is equivalent to say that the permutations that can be obtained in at most $p$ steps in the whole genome duplication - random loss model are exactly those whose number of descents is at most $2^p-1$.

Now, we can notice that the property of being obtainable in at most $p$ steps is stable for the pattern-involvement relation $\prec$: if $\sigma$ can be obtained in at most $p$ steps, and if $\pi \prec \sigma$, then $\pi$ can also be obtained in at most $p$ steps. Indeed, it is enough to perform the same duplication-loss scenario on $\sigma$, keeping track only of the elements of $\sigma$ that form an occurrence of $\pi$. This stability for $\prec$ implies that the class of permutations obtainable in at most $p$ steps is a class of pattern-avoiding permutations, whose excluded patterns are the minimal (again in the sense of $\prec$) permutations that cannot be obtained in $p$ steps.

Then, by Theorem \ref{thm:whole-genome}, the excluded patterns are the minimal permutations with $2^p$ descents. We initiated a study of the minimal permutations with $d$ descents in \cite{BP07}. However, it is simple to notice that a permutation with $d$ descents and minimal for this criterion has size at most $2d$, since it does not contain to consecutive ascents by minimality. An immediate consequence is that the number of excluded patterns is finite.

This allows us to state the following version of Theorem \ref{thm:whole-genome}:
\begin{thm}
The permutations that can be obtained in at most $p$ steps in the whole genome duplication - random loss model form a class of pattern-avoiding permutations. The excluded patterns are the permutations with exactly $2^p$ descents that are minimal (in the sense of $\prec$) for this criterion. These excluded patterns are in finite number.
\end{thm}

In \cite{BP07}, we will give a simpler description and some properties of these minimal permutations with $d$ descents.

\subsection{Permutations obtained in one step of width $K$}
\label{section:one-step}

As an introduction to the study of $\C(K,p)$, we deal in this section with the
simpler case of the class $\C(K) = \C(K,1)$ of permutations obtained from $12
\ldots n$ in one duplication-loss step of width at most $K$. Assume in this
section that the parameter $K \geq 2$ is fixed. Throughout this section, when
referring to a duplication-loss step, we always mean duplication-loss step of
width $K$, except when otherwise explicitly stated.

%
%

It is easily noticed that any permutation of $\C(K)$ cannot have more than one
descent. Conversly, any permutation of size at most $K$ having exactly one
descent belongs to $\C(K)$.

Although it is a technical point of importance in the proof of Theorem
\ref{thm:one-step}, the following proposition comes straightforward:
\begin{prop}
The permutations of size $K+1$ that do not belong to $\C(K)$ and having exactly one descent are exactly those of $S_{K+1}$ with one descent that do not start with $1$ nor end with $K+1$. \label{prop:one-step-sizeK+1}
\end{prop}

\begin{pf}
Let $\sigma = \sigma_1 \sigma_2 \ldots \sigma_{K+1}$ be a permutation of size $K+1$ that does not belong to $\C(K)$ but has exactly one descent. Now, if $\sigma_1 =1$, then $\overline{\sigma} = \sigma_2 \ldots \sigma_{K+1}$ is a permutation (of $\{2, 3, \ldots , K+1\}$) of size $K$ having one descent, and therefore $\overline{\sigma}$ can be obtained from $23\ldots K+1$ in one duplication-loss step. Applying the same transformation to $123\ldots K+1$ will then produce $\sigma$, contradicting that $\sigma \notin \C(K)$. The same reasoning holds when $\sigma_{K+1} = K+1$. So $\sigma$ does not start with $1$ nor end with $K+1$.

Now if $\sigma$ is a permutation of size $K+1$ having exactly one descent, that does not start with $1$ nor end with $K+1$, we claim that $\sigma$ cannot be obtained from $12 \ldots K+1$ in one duplication-loss step. This is because no duplication-loss step of width $K$ can move both $1$ and $K+1$ in $12 \ldots K+1$.
\end{pf}

\begin{thm}
The class $\C(K)$ of permutations obtained from $12 \ldots n$ (for some $n \geq 1$) in one duplication-loss step of width $K$ is a class $S(B)$ of pattern-avoiding permutations whose basis $B$ is composed of $3+2^{K-1}$ patterns of size at most $K+1$. Namely $B = \{321, 3142, 2143\} \cup D$, $D$ being the set of all permutations of $S_{K+1}$ that do not start with $1$ nor end with $K+1$, and containing exactly one descent.
\label{thm:one-step}
\end{thm}

\begin{exmp}
$ \C(4) =$ \\ $ S(321, 3142, 2143, 23451, 23514, 24513, 34512, 25134, 35124, 45123, 51234)$
\label{ex:thm-one-step}
\end{exmp}

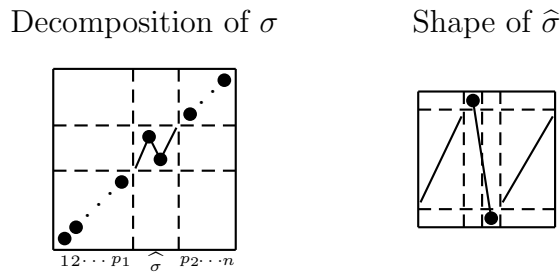
\begin{figure}[ht!]
\begin{center}
\psset{unit=0.3cm}
\begin{pspicture}(-4,0)(25,11)
\rput(4,10){Decomposition of $\sigma$}
\psline(0,0)(8,0)
\psline(0,0)(0,8)
\psline(8,8)(8,0)
\psline(8,8)(0,8)
\psline[linestyle=dashed](0,3.5)(8,3.5)
\psline[linestyle=dashed](0,5.5)(8,5.5)
\psline[linestyle=dashed](3.5,0)(3.5,8)
\psline[linestyle=dashed](5.5,0)(5.5,8)
\pscircle*(0.5,0.5){0.3}
\pscircle*(1,1){0.3}
\rput(1.5,1.5){.}
\rput(2,2){.}
\rput(2.5,2.5){.}
\pscircle*(3,3){0.3}
\rput(0.5,-0.5){{\tiny $1$}}
\rput(1,-0.5){{\tiny $2$}}
\rput(1.8,-0.5){{\tiny $\ldots$}}
\rput(3,-0.5){{\tiny $p_1$}}
\pscircle*(4.2,5){0.3}
\pscircle*(4.7,4){0.3}
\psline(3.6,3.6)(4.2,5)
\psline(4.2,5)(4.7,4)
\psline(4.7,4)(5.4,5.4)
\rput(4.5,-0.5){{\tiny $\widehat{\sigma}$}}
\pscircle*(6,6){0.3}
\rput(6.5,6.5){.}
\rput(7,7){.}
\pscircle*(7.5,7.5){0.3}
\rput(6,-0.5){{\tiny $p_2$}}
\rput(6.9,-0.5){{\tiny $\ldots$}}
\rput(7.7,-0.5){{\tiny $n$}}
\rput(19,10){Shape of $\widehat{\sigma}$}
\psline(16,1)(22,1)
\psline(16,1)(16,7)
\psline(22,7)(22,1)
\psline(22,7)(16,7)
\psline[linestyle=dashed](18,1)(18,7)
\psline[linestyle=dashed](18.8,1)(18.8,7)
\psline[linestyle=dashed](19.6,1)(19.6,7)
\psline[linestyle=dashed](16,1.8)(22,1.8)
\psline[linestyle=dashed](16,6.2)(22,6.2)
\pscircle*(19.2,1.4){0.3}
\pscircle*(18.4,6.6){0.3}
\psline(16.1,2.1)(17.9,5.9)
\psline(19.2,1.4)(18.4,6.6)
\psline(19.7,2.1)(21.9,5.9)
\end{pspicture}
\caption{Decomposition $\sigma = 12 \ldots p_1 \widehat{\sigma} p_2 (p_2 +1) \ldots n$ on the graphical representation of $\sigma$, and shape of $\widehat{\sigma}$ \label{fig:proof-thm-one-step}}
\end{center}
\end{figure}

\begin{pf}
We prove the reversed statement: $\sigma \notin S(B)$ if and only if $\sigma$ cannot be obtained from an identity permutation in one duplication-loss step of width $K$.

Assume $\sigma \notin S(B)$. Then there exists $b \in B$ such that $b \prec
\sigma$. If $b = 321$, $3142$ or $2143$, then $\sigma$ has at least $2$ descents and cannot be obtained in one duplication-loss step. Otherwise, using Proposition \ref{prop:one-step-sizeK+1}, there exists $\rho \in S_{K+1}$ such that $\rho \prec \sigma$ and $\rho \notin \C(K)$. Now if $\sigma$ could be obtained in one duplication-loss step, then so would be $\rho$, yielding a contradiction. So $\sigma \notin \C(K)$.

Conversly, assume that $\sigma \notin \C(K)$. If $\sigma$ contains at least $2$ descents, then $\sigma$ contains an occurrence of $321$ or $3142$ or $2143$, since these three are the minimal permutations (in the sense of the relation $\prec$) with $2$ descents. And consequently, $\sigma \notin S(B)$. Thus we may assume that $\sigma$ has exactly one descent. We decompose $\sigma \in S_n$ into $\sigma = 12 \ldots p_1 \widehat{\sigma} p_2 (p_2 +1) \ldots n$, where $\widehat{\sigma}$ is a permutation of the set $\{p_1+1 , p_1+2 \ldots, p_2-1\}$ that does not start with $p_1+1$ nor end with $p_2-1$, and contains exactly one descent. This decomposition is shown in Figure \ref{fig:proof-thm-one-step}. We denote by $\widehat{K}$ the size of $\widehat{\sigma}$. Since $\sigma \notin \C(K)$, necessarily $\widehat{K} \geq K+1$ or we would get a contradiction. If $\widehat{K} = K+1$, we get that $\widehat{\sigma}$ is an occurrence of some pattern of $D \subset B$ in $\sigma$. As a consequence, $\sigma \notin S(B)$. What is left to prove is that this extends to the case $\widehat{K} > K+1$. We just need to show that we can remove elements in $\widehat{\sigma}$ without violating any of the properties below:
\begin{itemize}
 \item the permutation does not start with its smallest element
 \item the permutation does not end with its greatest element
 \item the permutation has exactly one descent
\end{itemize}
until we get a permutation of size $K+1$. At that point $\widehat{\sigma}$ contains an occurrence of a pattern in $D$, and so does $\sigma$, and we get that $\sigma \notin S(B)$. Now, because of the conditions on $\widehat{\sigma}$, the only descent in $\widehat{\sigma}$ necessarily goes from the greatest to the smallest element in $\widehat{\sigma}$, ensuring that it is possible to remove elements without violating any of the properties above (see Figure \ref{fig:proof-thm-one-step}).
\end{pf}

\subsection{Permutations obtained in $p$ steps of width $K$}
\label{section:p-steps}

As for the case of $\C(K,1)$ in Section \ref{section:one-step}, we prove
(Theorem \ref{thm:p-steps}) in this section that the class $\C(K,p)$ of all permutations obtained from an identity permutation after $p$ duplication-loss steps of width at most $K$ is a class of pattern-avoiding permutations. However, we do not get a precise description of the basis of this class, but only an upper bound on the size of the excluded patterns. As in the previous section, when referring to a duplication-loss step, we always mean duplication-loss step of width $K$, except when otherwise explicitely stated.

To prove the announced result, we will need a few more notations and technical lemmas.

The \emph{vector} from $i$ to $j$ in a permutation $\sigma$ consists of all elements whose positions lie between the positions of $i$ and $j$, $i$ and $j$ being included. The \emph{size} of a vector is the number of elements in it. For example, the vector from $7$ to $2$ in the permutation $4123576$ is $\overleftarrow{2357}$, and has size $4$.

\begin{defn}
Let $\sigma$ be a permutation of $S_n$. The \emph{value-position vector associated with $i \in [1..n]$} (\emph{$vp$-vector} for short) is the vector of $\sigma$ going from $i$ to $\sigma_i$, if $i$ is not a fixpoint of $\sigma$. In the case $i = \sigma_i$, the $vp$-vector associated with $i$ is empty.
\label{def:vp-vector}
\end{defn}

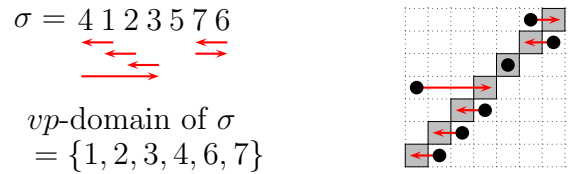
\begin{figure}[ht]
\begin{center}
\psset{unit=0.3cm}
\begin{pspicture}(-4,0)(23,10)
\rput(-1,7.5){$\sigma =$}
\rput(1,7.5){$4$}
\rput(2,7.5){$1$}
\rput(3,7.5){$2$}
\rput(4,7.5){$3$}
\rput(5,7.5){$5$}
\rput(6,7.5){$7$}
\rput(7,7.5){$6$}
\psline[linecolor=red]{->}(2.2,6.5)(0.8,6.5)
\psline[linecolor=red]{->}(3.2,6)(1.8,6)
\psline[linecolor=red]{->}(4.2,5.5)(2.8,5.5)
\psline[linecolor=red]{->}(0.8,5)(4.2,5)
\psline[linecolor=red]{->}(7.2,6.5)(5.8,6.5)
\psline[linecolor=red]{->}(5.8,6)(7.2,6)

\rput(3,3){$vp$-domain of $\sigma $}
\rput(4,1.5){$= \{1,2,3,4,6,7\}$ }
\psgrid[subgriddiv=1,gridwidth=.2pt,griddots=5,gridlabels=0pt](15,1)(22,8)
\psframe[fillstyle=solid,fillcolor=lightgray,linewidth=0](15,1)(16,2)
\psframe[fillstyle=solid,fillcolor=lightgray,linewidth=0](16,2)(17,3)
\psframe[fillstyle=solid,fillcolor=lightgray,linewidth=0](17,3)(18,4)
\psframe[fillstyle=solid,fillcolor=lightgray,linewidth=0](18,4)(19,5)
\psframe[fillstyle=solid,fillcolor=lightgray,linewidth=0](19,5)(20,6)
\psframe[fillstyle=solid,fillcolor=lightgray,linewidth=0](20,6)(21,7)
\psframe[fillstyle=solid,fillcolor=lightgray,linewidth=0](21,7)(22,8)
\psline[linecolor=red]{->}(16.5,1.5)(15.2,1.5)
\psline[linecolor=red]{->}(17.5,2.5)(16.2,2.5)
\psline[linecolor=red]{->}(18.5,3.5)(17.2,3.5)
\psline[linecolor=red]{->}(15.5,4.5)(18.8,4.5)
\psline[linecolor=red]{->}(21.5,6.5)(20.2,6.5)
\psline[linecolor=red]{->}(20.5,7.5)(21.8,7.5)
\pscircle*(15.5,4.5){0.3}
\pscircle*(16.5,1.5){0.3}
\pscircle*(17.5,2.5){0.3}
\pscircle*(18.5,3.5){0.3}
\pscircle*(19.5,5.5){0.3}
\pscircle*(20.5,7.5){0.3}
\pscircle*(21.5,6.5){0.3}
\end{pspicture}
\caption{$vp$-vectors and $vp$-domain for $\sigma = 4123576$, in the usual and in the graphical representations \label{fig:vp-vectors}}
\end{center}
\end{figure}

It should appear in this definition that the $vp$-vector associated with $i$, going from the element of $\sigma$ which has value $i$ to the element of $\sigma$ at position $i$, represents the necessary move for $i$ to reach its position in the sorted permutation $12\ldots n$. As it can be seen on Figure \ref{fig:vp-vectors}, on the graphical representation of permutations used throughout the paper, the $vp$-vector associated with $i$ is an arrow going horizontally from the element at ordinate $i$ to the diagonal.

We can also notice that a non-empty $vp$-vector contains at least two elements.

To take into account all the moves necessary to sort $\sigma$ to $12 \ldots n$, it is convenient to introduce the \emph{value-position domain}:
\begin{defn}
Let $\sigma$ be a permutation of $S_n$. The \emph{value-position domain} of $\sigma$ (\emph{$vp$-domain} for short) is composed of all elements of $\sigma$ appearing in at least one $vp$-vector.
\label{def:vp-domain}
\end{defn}

These two definitions are illustrated on Figure \ref{fig:vp-vectors}.

Now, observe that for any permutation, the $vp$-vectors are reversible in the sense that reversing all the arrows will give a set of vectors that represent the moves of elements that are necessary to "unsort" $12 \ldots n$ into $\sigma$. It is easily seen from Definitions \ref{def:vp-vector} and \ref{def:vp-domain} and this remark that for any permutation $\sigma \in \C(K,p)$, any element belonging to the $vp$-domain of $\sigma$ also belongs to at least one of the duplication-loss steps used to obtain $\sigma$ from $12 \ldots n$. Consequently, the $vp$-domain of $\sigma$ contains at most $Kp$ elements.

\begin{lem}
Consider a permutation $\sigma$, and the permutation $\tau$ obtained from $\sigma$ by the removal of some element $j$. Then for any element $i \neq j$ such that $i \neq \sigma_i$, either this element becomes a fixpoint in $\tau$ or the size of the $vp$-vector associated with this element in $\tau$ remains constant, is increased of 1 or is diminished of 1 with respect to the size of the $vp$-vector associated with $i$ in $\sigma$.
\label{lem:vp-vector-unchanged}
\end{lem}

\begin{figure}[ht]
\begin{center}
\psset{unit=0.28cm}
\begin{pspicture}(-2,-2)(40,12)
\psgrid[subgriddiv=1,gridwidth=.2pt,griddots=5,gridlabels=0pt](0,0)(8,8)
\psframe[fillstyle=solid,fillcolor=lightgray,linewidth=0](0,0)(1,1)
\psframe[fillstyle=solid,fillcolor=lightgray,linewidth=0](1,1)(2,2)
\psframe[fillstyle=solid,fillcolor=lightgray,linewidth=0](2,2)(3,3)
\psframe[fillstyle=solid,fillcolor=lightgray,linewidth=0](3,3)(4,4)
\psframe[fillstyle=solid,fillcolor=lightgray,linewidth=0](4,4)(5,5)
\psframe[fillstyle=solid,fillcolor=lightgray,linewidth=0](5,5)(6,6)
\psframe[fillstyle=solid,fillcolor=lightgray,linewidth=0](6,6)(7,7)
\psframe[fillstyle=solid,fillcolor=lightgray,linewidth=0](7,7)(8,8)
\psline[linecolor=blue,linewidth=0.1](2,-1)(2,9)
\psline[linecolor=blue,linewidth=0.1](3,-1)(3,9)
\psline[linecolor=blue,linewidth=0.1](-1,6)(9,6)
\psline[linecolor=blue,linewidth=0.1](-1,7)(9,7)
\pscircle[fillstyle=solid,fillcolor=blue,linecolor=blue](2.5,6.5){0.3}
\rput(-1.5,6.5){\color{blue}$j$}
\psframe[fillstyle=crosshatch,hatchcolor=blue,linecolor=blue,hatchsep=0.2](3,2)(4,3)
\psframe[fillstyle=crosshatch,hatchcolor=blue,linecolor=blue,hatchsep=0.2](4,3)(5,4)
\psframe[fillstyle=crosshatch,hatchcolor=blue,linecolor=blue,hatchsep=0.2](5,4)(6,5)
\psframe[fillstyle=crosshatch,hatchcolor=blue,linecolor=blue,hatchsep=0.2](6,5)(7,6)
\psframe[fillstyle=solid,fillcolor=black,linewidth=0](3,4)(4,5)
\psframe[fillstyle=solid,fillcolor=black,linewidth=0](3,5)(4,6)
\psframe[fillstyle=solid,fillcolor=black,linewidth=0](4,5)(5,6)
\psframe[fillstyle=vlines,hatchcolor=gray,linewidth=0.08,hatchangle=-45,hatchsep=0.2](0,7)(1,8)
\psframe[fillstyle=vlines,hatchcolor=gray,linewidth=0.08,hatchangle=-45,hatchsep=0.2](1,7)(2,8)
\psframe[fillstyle=vlines,hatchcolor=gray,linewidth=0.08,hatchangle=-45,hatchsep=0.2](3,0)(4,1)
\psframe[fillstyle=vlines,hatchcolor=gray,linewidth=0.08,hatchangle=-45,hatchsep=0.2](4,0)(5,1)
\psframe[fillstyle=vlines,hatchcolor=gray,linewidth=0.08,hatchangle=-45,hatchsep=0.2](5,0)(6,1)
\psframe[fillstyle=vlines,hatchcolor=gray,linewidth=0.08,hatchangle=-45,hatchsep=0.2](6,0)(7,1)
\psframe[fillstyle=vlines,hatchcolor=gray,linewidth=0.08,hatchangle=-45,hatchsep=0.2](7,0)(8,1)
\psframe[fillstyle=vlines,hatchcolor=gray,linewidth=0.08,hatchangle=-45,hatchsep=0.2](3,1)(4,2)
\psframe[fillstyle=vlines,hatchcolor=gray,linewidth=0.08,hatchangle=-45,hatchsep=0.2](4,1)(5,2)
\psframe[fillstyle=vlines,hatchcolor=gray,linewidth=0.08,hatchangle=-45,hatchsep=0.2](5,1)(6,2)
\psframe[fillstyle=vlines,hatchcolor=gray,linewidth=0.08,hatchangle=-45,hatchsep=0.2](6,1)(7,2)
\psframe[fillstyle=vlines,hatchcolor=gray,linewidth=0.08,hatchangle=-45,hatchsep=0.2](7,1)(8,2)
\psframe[fillstyle=vlines,hatchcolor=gray,linewidth=0.08,hatchangle=-45,hatchsep=0.2](4,2)(5,3)
\psframe[fillstyle=vlines,hatchcolor=gray,linewidth=0.08,hatchangle=-45,hatchsep=0.2](5,2)(6,3)
\psframe[fillstyle=vlines,hatchcolor=gray,linewidth=0.08,hatchangle=-45,hatchsep=0.2](6,2)(7,3)
\psframe[fillstyle=vlines,hatchcolor=gray,linewidth=0.08,hatchangle=-45,hatchsep=0.2](7,2)(8,3)
\psframe[fillstyle=vlines,hatchcolor=gray,linewidth=0.08,hatchangle=-45,hatchsep=0.2](5,3)(6,4)
\psframe[fillstyle=vlines,hatchcolor=gray,linewidth=0.08,hatchangle=-45,hatchsep=0.2](6,3)(7,4)
\psframe[fillstyle=vlines,hatchcolor=gray,linewidth=0.08,hatchangle=-45,hatchsep=0.2](7,3)(8,4)
\psframe[fillstyle=vlines,hatchcolor=gray,linewidth=0.08,hatchangle=-45,hatchsep=0.2](6,4)(7,5)
\psframe[fillstyle=vlines,hatchcolor=gray,linewidth=0.08,hatchangle=-45,hatchsep=0.2](7,4)(8,5)
\psframe[fillstyle=vlines,hatchcolor=gray,linewidth=0.08,hatchangle=-45,hatchsep=0.2](7,5)(8,6)
\psline[linecolor=red,linewidth=0.15]{->}(6.5,2.5)(5.5,2.5)
\psline[linecolor=red,linewidth=0.15]{->}(4.5,5.5)(3.5,5.5)
\psline[linecolor=red,linewidth=0.15]{->}(0.5,7.5)(0.5,6.5)
\psline[linecolor=red,linewidth=0.15]{->}(5.5,7.5)(4.5,6.5)
\pscircle[fillstyle=solid,fillcolor=red,linecolor=red](1.5,0.5){0.15}
\pscircle[fillstyle=solid,fillcolor=red,linecolor=red](0.5,3.5){0.15}
\psgrid[subgriddiv=1,gridwidth=.2pt,griddots=5,gridlabels=0pt](14,0)(22,8)
\psframe[fillstyle=solid,fillcolor=lightgray,linewidth=0](14,0)(15,1)
\psframe[fillstyle=solid,fillcolor=lightgray,linewidth=0](15,1)(16,2)
\psframe[fillstyle=solid,fillcolor=lightgray,linewidth=0](16,2)(17,3)
\psframe[fillstyle=solid,fillcolor=lightgray,linewidth=0](17,3)(18,4)
\psframe[fillstyle=solid,fillcolor=lightgray,linewidth=0](18,4)(19,5)
\psframe[fillstyle=solid,fillcolor=lightgray,linewidth=0](19,5)(20,6)
\psframe[fillstyle=solid,fillcolor=lightgray,linewidth=0](20,6)(21,7)
\psframe[fillstyle=solid,fillcolor=lightgray,linewidth=0](21,7)(22,8)
\psline[linecolor=blue,linewidth=0.1](19,-1)(19,9)
\psline[linecolor=blue,linewidth=0.1](20,-1)(20,9)
\psline[linecolor=blue,linewidth=0.1](13,2)(23,2)
\psline[linecolor=blue,linewidth=0.1](13,3)(23,3)
\pscircle[fillstyle=solid,fillcolor=blue,linecolor=blue](19.5,2.5){0.3}
\rput(12.5,2.5){\color{blue}$j$}
\psframe[fillstyle=crosshatch,hatchcolor=blue,linecolor=blue,hatchsep=0.2](16,3)(17,4)
\psframe[fillstyle=crosshatch,hatchcolor=blue,linecolor=blue,hatchsep=0.2](17,4)(18,5)
\psframe[fillstyle=crosshatch,hatchcolor=blue,linecolor=blue,hatchsep=0.2](18,5)(19,6)
\psframe[fillstyle=solid,fillcolor=black,linewidth=0](18,3)(19,4)
\psframe[fillstyle=vlines,hatchcolor=gray,linewidth=0.08,hatchangle=-45,hatchsep=0.2](20,0)(21,1)
\psframe[fillstyle=vlines,hatchcolor=gray,linewidth=0.08,hatchangle=-45,hatchsep=0.2](21,0)(22,1)
\psframe[fillstyle=vlines,hatchcolor=gray,linewidth=0.08,hatchangle=-45,hatchsep=0.2](20,1)(21,2)
\psframe[fillstyle=vlines,hatchcolor=gray,linewidth=0.08,hatchangle=-45,hatchsep=0.2](21,1)(22,2)
\psframe[fillstyle=vlines,hatchcolor=gray,linewidth=0.08,hatchangle=-45,hatchsep=0.2](14,7)(15,8)
\psframe[fillstyle=vlines,hatchcolor=gray,linewidth=0.08,hatchangle=-45,hatchsep=0.2](15,7)(16,8)
\psframe[fillstyle=vlines,hatchcolor=gray,linewidth=0.08,hatchangle=-45,hatchsep=0.2](16,7)(17,8)
\psframe[fillstyle=vlines,hatchcolor=gray,linewidth=0.08,hatchangle=-45,hatchsep=0.2](17,7)(18,8)
\psframe[fillstyle=vlines,hatchcolor=gray,linewidth=0.08,hatchangle=-45,hatchsep=0.2](18,7)(19,8)
\psframe[fillstyle=vlines,hatchcolor=gray,linewidth=0.08,hatchangle=-45,hatchsep=0.2](14,6)(15,7)
\psframe[fillstyle=vlines,hatchcolor=gray,linewidth=0.08,hatchangle=-45,hatchsep=0.2](15,6)(16,7)
\psframe[fillstyle=vlines,hatchcolor=gray,linewidth=0.08,hatchangle=-45,hatchsep=0.2](16,6)(17,7)
\psframe[fillstyle=vlines,hatchcolor=gray,linewidth=0.08,hatchangle=-45,hatchsep=0.2](17,6)(18,7)
\psframe[fillstyle=vlines,hatchcolor=gray,linewidth=0.08,hatchangle=-45,hatchsep=0.2](18,6)(19,7)
\psframe[fillstyle=vlines,hatchcolor=gray,linewidth=0.08,hatchangle=-45,hatchsep=0.2](14,5)(15,6)
\psframe[fillstyle=vlines,hatchcolor=gray,linewidth=0.08,hatchangle=-45,hatchsep=0.2](15,5)(16,6)
\psframe[fillstyle=vlines,hatchcolor=gray,linewidth=0.08,hatchangle=-45,hatchsep=0.2](16,5)(17,6)
\psframe[fillstyle=vlines,hatchcolor=gray,linewidth=0.08,hatchangle=-45,hatchsep=0.2](17,5)(18,6)
\psframe[fillstyle=vlines,hatchcolor=gray,linewidth=0.08,hatchangle=-45,hatchsep=0.2](14,4)(15,5)
\psframe[fillstyle=vlines,hatchcolor=gray,linewidth=0.08,hatchangle=-45,hatchsep=0.2](15,4)(16,5)
\psframe[fillstyle=vlines,hatchcolor=gray,linewidth=0.08,hatchangle=-45,hatchsep=0.2](16,4)(17,5)
\psframe[fillstyle=vlines,hatchcolor=gray,linewidth=0.08,hatchangle=-45,hatchsep=0.2](14,3)(15,4)
\psframe[fillstyle=vlines,hatchcolor=gray,linewidth=0.08,hatchangle=-45,hatchsep=0.2](15,3)(16,4)
\psline[linecolor=red,linewidth=0.15]{->}(21.5,0.5)(20.5,0.5)
\psline[linecolor=red,linewidth=0.15]{->}(18.5,3.5)(18.5,2.5)
\psline[linecolor=red,linewidth=0.15]{->}(15.5,6.5)(15.5,5.5)
\psline[linecolor=red,linewidth=0.15]{->}(21.5,5.5)(20.5,4.5)
\psline[linecolor=red,linewidth=0.15]{->}(20.5,7.5)(19.5,6.5)
\pscircle[fillstyle=solid,fillcolor=red,linecolor=red](14.5,1.5){0.15}
\pscircle[fillstyle=solid,fillcolor=red,linecolor=red](16.5,0.5){0.15}
\psframe[fillstyle=solid,fillcolor=lightgray,linewidth=0](26,8)(27,9)
\rput(30.9,8.3){Diagonal}
\psframe[fillstyle=crosshatch,hatchcolor=blue,linecolor=blue,hatchsep=0.2](26,6)(27,7)
\rput(33.7,6.3){Candidate fixpoints}
\psline[linecolor=red,linewidth=0.15]{->}(26,5)(27,5)
\pscircle[fillstyle=solid,fillcolor=red,linecolor=red](26.5,4.3){0.15}
\rput(30.8,4.8){Changes}
\rput(33,3.8){in the $vp$-vectors}
\rput(33,2.1){Variation of the distance}
\rput(33,1.1){to the diagonal}
\psframe[fillstyle=solid,fillcolor=white,linewidth=0](26,-0.5)(27,0.5)
\rput(28.2,-0.2){$0$}
\psframe[fillstyle=solid,fillcolor=black,linewidth=0](29.7,-0.5)(30.7,0.5)
\rput(32.2,-0.2){$+1$}
\psframe[fillstyle=vlines,hatchcolor=gray,linewidth=0.08,hatchangle=-45,hatchsep=0.2](33.7,-0.5)(34.7,0.5)
\rput(36.2,-0.2){$-1$}
\end{pspicture}
\caption{Variation of the size of $vp$-vectors due to the removal of an element $j$ above or below the diagonal. \label{fig:vp-vector-unchanged}}
\end{center}
\end{figure}
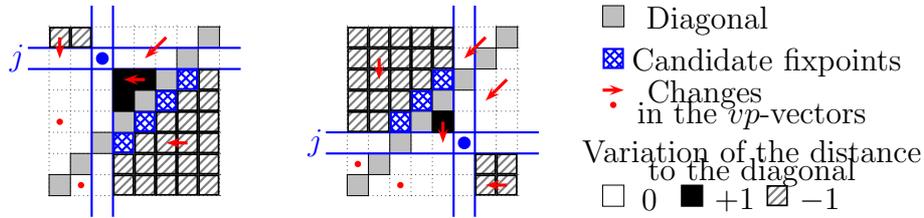

\begin{pf}
It is easily seen on the graphical representation of $\sigma$. Any element that does not lie just above or just below the diagonal cannot become a fixpoint when removing an element $j$. For elements that do not becom fixpoints, the horizontal distance to the diagonal can only change of $0$, $1$ or $-1$ when removing some element $j$ (see Figure \ref{fig:vp-vector-unchanged}).
\end{pf}

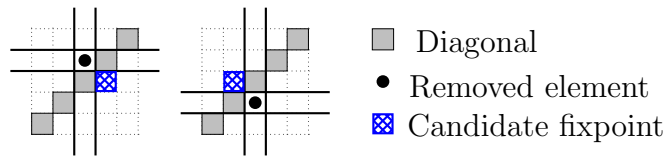
\begin{figure}[ht]
\begin{center}
\psset{unit=0.28cm}
\begin{pspicture}(-2,-2)(30,7)
\psgrid[subgriddiv=1,gridwidth=.2pt,griddots=5,gridlabels=0pt](0,0)(5,5)
\psframe[fillstyle=solid,fillcolor=lightgray,linewidth=0](0,0)(1,1)
\psframe[fillstyle=solid,fillcolor=lightgray,linewidth=0](1,1)(2,2)
\psframe[fillstyle=solid,fillcolor=lightgray,linewidth=0](2,2)(3,3)
\psframe[fillstyle=solid,fillcolor=lightgray,linewidth=0](3,3)(4,4)
\psframe[fillstyle=solid,fillcolor=lightgray,linewidth=0](4,4)(5,5)
\pscircle*(2.5,3.5){0.3}
\psline[linewidth=0.1](2,-1)(2,6)
\psline[linewidth=0.1](3,-1)(3,6)
\psline[linewidth=0.1](-1,3)(6,3)
\psline[linewidth=0.1](-1,4)(6,4)
\psframe[fillstyle=crosshatch,hatchcolor=blue,linecolor=blue,hatchsep=0.2](3,2)(4,3)
\psgrid[subgriddiv=1,gridwidth=.2pt,griddots=5,gridlabels=0pt](8,0)(13,5)
\psframe[fillstyle=solid,fillcolor=lightgray,linewidth=0](8,0)(9,1)
\psframe[fillstyle=solid,fillcolor=lightgray,linewidth=0](9,1)(10,2)
\psframe[fillstyle=solid,fillcolor=lightgray,linewidth=0](10,2)(11,3)
\psframe[fillstyle=solid,fillcolor=lightgray,linewidth=0](11,3)(12,4)
\psframe[fillstyle=solid,fillcolor=lightgray,linewidth=0](12,4)(13,5)
\pscircle*(10.5,1.5){0.3}
\psline[linewidth=0.1](10,-1)(10,6)
\psline[linewidth=0.1](11,-1)(11,6)
\psline[linewidth=0.1](7,1)(14,1)
\psline[linewidth=0.1](7,2)(14,2)
\psframe[fillstyle=crosshatch,hatchcolor=blue,linecolor=blue,hatchsep=0.2](9,2)(10,3)
\psframe[fillstyle=solid,fillcolor=lightgray,linewidth=0](16,4)(17,5)
\rput(20.9,4.3){Diagonal}
\pscircle*(16.5,2.5){0.3}
\rput(23.4,2.3){Removed element}
\psframe[fillstyle=crosshatch,hatchcolor=blue,linecolor=blue,hatchsep=0.2](16,0)(17,1)
\rput(23.7,0.3){Candidate fixpoint}
\end{pspicture}
\caption{The only fixpoint that can appear when removing a quasi-diagonal element. \label{fig:one-fixpoint}}
\end{center}
\end{figure}

\begin{lem}
For any permutation $\sigma$, there is at least one element $j$ such that the permutation $\tau$ obtained from $\sigma$ by the removal of $j$ contains at most one more fixpoint than $\sigma$.
\label{lem:one-fixpoint}
\end{lem}
\begin{pf}
It is convenient to introduce the \emph{quasi-diagonal elements} of $\sigma$, defined as follows. $i$ is a quasi-diagonal element of $\sigma$ if $\sigma_{i-1} = i$ or $\sigma_{i+1} =i$. These two cases correspond respectively to elements of $\sigma$ lying just above or just below the diagonal in the graphical representation of $\sigma$. Any element of $\sigma$ that may become a fixpoint in $\tau$ is necessarily a quasi-diagonal element.

If there is no quasi-diagonal element, then we can remove any element $j$ to obtain a permutation $\tau$ that does not have more fixpoints than $\sigma$. If there are some, then we pick $j$ among the quasi-diagonal elements. We claim that at most one fixpoint is create while removing $j$. The argument is simple. Suppose $j$ is such that $\sigma_{j-1} = j$, the other case being similar. Then the only fixpoint that may appear is $j-1$, if $\sigma_j = j-1$. This should appear clearly on Figure \ref{fig:one-fixpoint}.
\end{pf}

\begin{lem}
Consider a permutation $\sigma \notin \C(K,p)$ such that for any strict pattern $\tau$ of $\sigma$, $\tau \in \C(K,p)$. Then the $vp$-domain of $\sigma$ is of size at most $2Kp+2$.
\label{lem:bounded-vp-domain}
\end{lem}
\begin{pf}
By Lemma \ref{lem:one-fixpoint}, we can choose some $\tau \prec \sigma$ with $|\tau|+1 = |\sigma|$ and such that $\tau$ has at most one more fixpoint than $\sigma$. Call $j$ the element deleted in $\sigma$ to obtain $\tau$. By a previous remark, since $\tau \in \C(K,p)$, the $vp$-domain of $\tau$ is of size at most $Kp$, and is therefore composed of at most $Kp$ $vp$-vectors. Each of these $vp$-vectors in $\tau$ yields a $vp$-vector in $\sigma$, whose size is smaller or equal or possibly increased by $1$. Let us denote by $\mathcal{\overrightarrow{V}}$ the set of $vp$-vectors of $\sigma$ obtained from a $vp$-vector of $\tau$. Then the number of elements of $\sigma$ that belong to a $vp$-vector of $\mathcal{\overrightarrow{V}}$ is at most $2Kp$. However $\mathcal{\overrightarrow{V}}$ is not yet the $vp$-domain of $\sigma$. We must complete it with up to two $vp$-vectors: the one associated with the element $j$ deleted, and the one associated with the fixpoint of $\tau$ that was not a fixpoint in $\sigma$, if such a point exists. If such an element exists, then it is a quasi-diagonal element in $\sigma$ and its $vp$-vector (denoted $\overrightarrow{v}$) in $\sigma$ is necessarily of size $2$, so that $\mathcal{\overrightarrow{V}} \cup \{\overrightarrow{v}\}$ has total size at most $2Kp+2$. Now it is easily observed that any element of $\sigma$ belonging to one $vp$-vector necessarily belongs to at least two $vp$-vectors (this can be seen as a ``balance condition''). Consequently, all the elements of the $vp$-vector associated with $j$ are already covered by a vector of $\mathcal{\overrightarrow{V}} \cup \{\overrightarrow{v}\}$, so that the $vp$-domain of $\sigma$ is exactly the set of elements covered by $\mathcal{\overrightarrow{V}} \cup \{\overrightarrow{v}\}$. Therefore, its size is at most $2Kp+2$.
\end{pf}

\begin{figure}[ht]
\begin{center}
\psset{unit=0.3cm}
\begin{pspicture}(-1.2,-3)(16,3)
\rput(-1.2,0.3){$\sigma =$}
\psframe[fillstyle=solid,fillcolor=red,linewidth=0](1.5,0)(2,1)
\psframe[fillstyle=solid,fillcolor=red,linewidth=0](4,0)(7,1)
\psframe[fillstyle=solid,fillcolor=red,linewidth=0](8.5,0)(9.5,1)
\psframe[fillstyle=solid,fillcolor=red,linewidth=0](12.5,0)(14.5,1)
\psframe(0,0)(15,1)
\psline[linecolor=red]{->}(1.75,1.1)(7.35,2.3)
\psline[linecolor=red]{->}(5.5,1.1)(7.45,2.3)
\psline[linecolor=red]{->}(9,1.1)(7.55,2.3)
\psline[linecolor=red]{->}(13.5,1.1)(7.65,2.3)
\rput(7.5,3){at most $Kp+1$ $vp$-windows}
\psline{->}(0.75,-0.1)(7.3,-2.3)
\psline{->}(3,-0.1)(7.4,-2.3)
\psline{->}(7.75,-0.1)(7.5,-2.3)
\psline{->}(11,-0.1)(7.6,-2.3)
\psline{->}(14.75,-0.1)(7.7,-2.3)
\rput(7.5,-3){at most $Kp+2$ free-windows}
\end{pspicture}
\caption{Proof of Lemma \ref{lem:nearly-sorted-pi} \label{fig:nearly-sorted-pi}}
\end{center}
\end{figure}
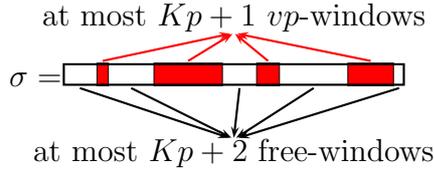

\begin{lem}
Consider a permutation $\sigma \notin \C(K,p)$ of size $n > (Kp+2)^2-2$ such that for any strict pattern $\tau$ of $\sigma$, $\tau \in \C(K,p)$. Then $\sigma$ is of the form $\sigma = I i (i+1) \ldots (i+Kp) J$ with $I$ a permutation of $[1..i-1]$ and $J$ a permutation of $[i+Kp+1..n]$. It is possible that $I$ or $J$ is empty.
\label{lem:nearly-sorted-pi}
\end{lem}
\begin{pf}
By Lemma \ref{lem:bounded-vp-domain}, the $vp$-domain of $\sigma$ is of size at
most $2Kp+2$. We can decompose $\sigma$ into \emph{free windows} of consecutive elements outside the $vp$-domain of $\sigma$, separated by windows of consecutive elements of the $vp$-domain. Now, there are at most $Kp+1$ windows of consecutive elements of the $vp$-domain, and consequently, there are at most $Kp+2$ free windows in $\sigma$. Since $\sigma$ is of size $n > (Kp+2)^2-2 = (Kp+2)Kp + 2Kp+2$, at least one of the free windows of $\sigma$ has size strictly greater than $Kp$, \emph{i.e.,} contains at least $Kp+1$ elements. By definition, these elements do not belong to the $vp$-domain of $\sigma$, and hence they allow the decomposition of $\sigma$ into $\sigma = I i (i+1) \ldots (i+Kp) J$ with $I$ a permutation of $[1..i-1]$ and $J$ a permutation of $[i+Kp+1..n]$. Figure \ref{fig:nearly-sorted-pi} represent the decomposition of $\sigma$ used in this proof.
\end{pf}

\begin{lem}
Consider a permutation $\sigma = \sigma'(j+1)(j+2)\ldots n$ where $\sigma'$ is
a permutation of $[1\ldots j]$. If $\sigma$ is obtainable after $p$ duplication
steps of size at most $K$ then $\sigma$ is obtainable after $p$ duplication steps of size at
most $K$ such that the duplicated window for each step does not intersect 
\label{lem:localisation}
\end{lem}

\begin{pf}
The key idea is to consider the first sequence $s_1,s_2,\ldots,s_p$ of
duplication-loss steps and create a new sequence $s'_1,s'_2,\ldots,s'_p$ such
that :
\begin{itemize}
  \item Each step $s'_i$ concerns only elements of $[1..j]$.
  \item After every step $s'_i$, the elements $1,2,\ldots,j$ are in the same
  order than after performing steps $s_1,s_2,\ldots,s_i$.
\end{itemize}
Then the proof is by induction on the number of steps. If there is only one step
then the proof is straighforward. Suppose now that the above statement is true
until $p-1$ steps. Then for the last step, we use our hypothesis for $p-1$ so
that we have operations $s'_1,s'_2,\ldots,s'_{p-1}$ respecting the above
conditions. For $s'_n$, only notice that the elements of $[1\ldots j]$ involved
in $s_n$ are also in a window of size $K$ in the permutation obtained after
$s'_{j-1}$ and in the same relative order by our induction hypothesis which
proves the existence of $s'_n$.
\end{pf}

Using these lemmas, we state and prove a key proposition:

\begin{prop}
Consider a permutation $\sigma \notin \C(K,p)$. Then either $\sigma$ is of size at most $(Kp+2)^2-2$, or there exists a strict pattern $\tau$ of $\sigma$ that does not belong to $\C(K,p)$.
\label{prop:p-steps}
\end{prop}
\begin{pf}
Consider a permutation $\sigma \notin \C(K,p)$ such that any strict pattern $\tau$ of $\sigma$ belongs to $\C(K,p)$. We want to show that $\sigma$ is of size $n \leq (Kp+2)^2-2$. Let us assume the contrary. By Lemma \ref{lem:nearly-sorted-pi}, there exist $i \in [1..n]$, $I$ a permutation of $[1..i-1]$ and $J$ a permutation of $[i+Kp+1..n]$ such that $\sigma = I i (i+1) \ldots (i+Kp) J$. Let us denote $\widehat{\sigma}$ the permutation $\widehat{\sigma} = I i (i+1) \ldots (i+Kp-1) (J-1)$, where $(J-1)$ is the permutation of $[i+Kp..n-1]$ obtained from $J$ by subtracting $1$ to every element of $J$. $\widehat{\sigma}$ is a strict pattern of $\sigma$, hence $\widehat{\sigma} \in \C(K,p)$. Consider a shortest sequence of duplication-loss steps of width at most $K$ that produces $\widehat{\sigma}$ from $12\ldots(n-1)$. This sequence has at most $p$ steps, each of width at most $K$. It implies that the total distance crossed by the elements that are duplicated is at most $Kp$. Consequently, it is not possible to bring an element of $I$ and an element of $J-1$ in two consecutive positions. So it is necessary that the duplication-loss steps of the scenario we consider are \emph{internal} to $I$ and $J-1$. We can reproduce these steps in $I$ and $J$ to obtain $\sigma$ from $12\ldots n$ in at most $p$ duplication-loss steps of width at most $K$, contradicting that $\sigma \notin \C(K,p)$.
\end{pf}

It is then quite easy to prove Theorem \ref{thm:p-steps}: 

\begin{thm}
The class $\C(K,p)$ of all permutations obtained from an identity permutation after $p$ duplication-loss steps of width at most $K$ is a class of pattern-avoiding permutations whose basis is finite and contains only patterns of size at most $(Kp+2)^2-2$.
\label{thm:p-steps}
\end{thm}
\begin{pf}
We set $B = \{\pi : \pi \notin \C(K,p) \textrm{ and } |\pi| \leq (Kp+2)^2-2 \}$ and show that $S(B) = \C(K,p)$.

Consider $\sigma \notin \C(K,p)$. If $|\sigma| \leq (Kp+2)^2-2$, then $\sigma \in
B$ and $\sigma \notin S(B)$. Otherwise, if $|\sigma| > (Kp+2)^2-2$, then by
Proposition \ref{prop:p-steps}, there exists a strict pattern $\tau$ of $\sigma$
that does not belong to $\C(K,p)$. Reasoning by induction on the size of the
permutations, we deduce from $\tau \notin \C(K,p)$ that $\tau \notin S(B)$. A
direct consequence is that $\sigma \notin S(B)$. This proves that $S(B) \subseteq
\C(K,p)$.

Conversely, consider $\sigma \in \C(K,p)$. Then any pattern $\tau$ of $\sigma$
is also obtainable from an identity permutation in at most $p$ steps of width at most $K$ (using the sequence of duplication-loss steps associated with $\sigma$), \emph{i.e.,} $\tau \in \C(K,p)$. Then $\sigma$ does not contain an occurrence of any pattern of $B$, \emph{i.e.,} $\sigma \in  S(B)$. This shows that $\C(K,p) \subseteq S(B)$, ending the proof of the theorem.
\end{pf}

\section{Number of steps of width $K$ to obtain any permutation of size $n$}
\label{section:algo}

The \emph{whole genome duplication - random loss} model is studied in
\cite{CCMR06}, and the authors describe a method to compute an optimal
duplication-loss scenario, \emph{i.e.,} a scenario of duplications (of the
whole genome in this case) and losses whose number of steps is minimal.

Our model with bounded size duplication operations reduces to the \emph{whole
genome duplication - random loss} case when $K = n$ and thus to a radix-sort
algorithm as shown in \cite{CCMR06} and to a bubble-sort when $K=2$. Thus we give
some algorithm whose complexity matches the two extremal cases and shows some
continuity between the two sorting algorithms.

It is worth noticing that any scenario in our model can be viewed as a whole
genome duplication - random loss scenario. Consequently, the number of steps of
an optimal whole genome duplication - random loss scenario is a lower bound to
the number of steps of an optimal scenario in our duplication-loss model.

It is also easy to see that, when considering permutations of size at most $K$,
our model and the whole genome duplication - random loss model coincide. Indeed,
we will use for our purpose the procedure of \cite{CCMR06}, which is given in
Algorithm \ref{alg:binary}. We omit the proof of correctness and optimality of
this algorithm. See \cite{CCMR06} for details.

\begin{algorithm}[ht]
\caption{An optimal whole genome duplication - random loss scenario from $12 \ldots K$ to $\sigma \in S_K$}
\label{alg:binary}
\begin{algorithmic}[1]
\STATE {$\pi = 12 \ldots K$}
\STATE {Partition $\sigma$ into maximal increasing substrings, from left to right}
\STATE {Each element of $[1..K]$ appearing in the $i^{th}$ maximal increasing substring gets as a label the binary representation of $i$}
\FOR{$j = 1$ to $\lceil \log_2(desc(\sigma)+1) \rceil$}
\STATE {Perform a duplication-loss step on $\pi$ that keeps in the first copy of $\pi$ exactly the elements whose label has a $0$ in its $j^{th}$ least significant bit}
\ENDFOR
\end{algorithmic}
\end{algorithm}

In order to examine every bit of the labels given to the elements of $[1..K]$, the number of steps in the loop on line 4 is $\lceil \log_2($number of maximal increasing substrings of $\sigma) \rceil = \lceil \log_2(desc(\sigma)+1) \rceil$. A consequence is that the number of steps in an optimal whole genome duplication - random loss scenario from $12 \ldots n$ to $\sigma$ is $\Theta(\log n)$ in the worst case and on average (see equation (\ref{eq:log-average}) for the average case).

Note that the same algorithm can be used to compute an optimal whole genome duplication - random loss scenario from $i_1 i_2 \ldots i_k$, with $k \leq K$ and $i_1 < i_2 < \ldots < i_k$, to any permutation of $\{i_1, i_2, \ldots, i_k\}$.

\subsection{Upper bound}

In this section, we provide an algorithm that computes, for any permutation $\sigma \in S_n$ in input, a possible scenario of duplications and losses to obtain $\sigma$ from $12 \ldots n$. We will restrict ourselves to duplication-loss steps of width at most $K$, so that the number of duplication-loss steps corresponds to the cost of the scenario in our cost model. We are interested in the number of duplication-loss steps of the scenario produced by the algorithm, in the worst case, and on average. It provides an upper bound on the number of duplication-loss steps that are necessary to obtain a permutation. The algorithm we use is described in Algorithm \ref{alg:bucket}.

\begin{algorithm}[ht]
\caption{A duplication-loss scenario from $12 \ldots n$ to $\sigma \in S_n$}
\label{alg:bucket}
\begin{algorithmic}[1]
\STATE {$\pi \leftarrow 12 \ldots n$}
\FOR{$i = 1$ to $\lceil \frac{n-K}{\lfloor K/2 \rfloor} \rceil$ }
\STATE {Let $L^i = \{ \sigma_{j} : n - i \lfloor K/2 \rfloor +1 \leq j \leq n-(i-1) \lfloor K/2 \rfloor \}$}
\STATE {Perform duplication-loss steps on $\pi$ to move from left to right the elements of $L^i$ to the positions $n - i \lfloor K/2 \rfloor +1$ to $n-(i-1) \lfloor K/2 \rfloor$ of $\pi$, without changing their respective order}
\ENDFOR
\FOR{$i = 1$ to $\lceil \frac{n-K}{\lfloor K/2 \rfloor} \rceil$}
\STATE {Perform Algorithm \ref{alg:binary} on the window of $\pi$ between the indices $n-i \lfloor K/2 \rfloor +1$ and $n-(i-1)\lfloor K/2 \rfloor$}
\ENDFOR
\STATE {Perform Algorithm \ref{alg:binary} on the window of $\pi$ between the indices $1$ and $n- \lceil \frac{n-K}{\lfloor K/2 \rfloor} \rceil \lfloor K/2 \rfloor$}
\end{algorithmic}
\end{algorithm}

A few keys to understand Algorithm \ref{alg:bucket} are the following remarks.

The set $L^i$ of values defined at line 3 represents the rightmost $\lfloor
K/2 \rfloor$ elements of $\sigma$ not yet examined. 
The algorithm consists in two different loops, the first one corresponding to
lines $2$ to $5$ of the algorithm and the second one from line $6$ to $8$.
At the end of the first loop (line 5), $\pi$ is decomposed into windows of width $\lfloor K/2 \rfloor$ (except the leftmost one which is of width at most $K$) ; and each of these windows is an increasing sequence containing exactly the same elements as the window of $\sigma$ corresponding to the same indices. In the second loop, we consider these windows from right to left and since there are of width less than $K$, we can call Algorithm \ref{alg:binary} (that implements whole genome duplication-random loss) on each window successively to transform $\pi$ into $\sigma$.

An example is given with $\sigma = 2\, 10\, 1\, 7\, 6\, 5\, 8\, 9\, 3\,
4$ and $K = 6$. 
We first cut $\sigma$ in chunks of size $\lfloor K/2 \rfloor = 3$ and obtain
$2\, 10\, 1\, |\, 7\, 6\, 5\, | \, 8\, 9\, 3\, | \,4$.
Then the first loop of the algorithm (step $2$ to $5$) starts from $1\, 2\,
3\, 4\, 5\, 6\, 7\, 8\, 9\, 10$ and takes the elements in increasing order to the
same chunk the belong to in $\sigma$. This gives $1\, 2\, 10\, |\,
5\, 6\, 7\, | \, 3\, 8\, 9\, | \,4$.
Then the second loop sorts each chunk separately to obtain $\sigma$ using the
radix sort Algorithm \ref{alg:binary} introduced in \cite{CCMR06}.

Notice here that we use in the second loop (except for the leftmost window) only
duplication-loss steps of width $\lfloor K/2 \rfloor$. An improvement we
considered is to use whole genome duplication - random loss scenarios on windows
of width $K$, that are nonetheless increasing sequences. Unfortunately, we were
not able to analyse how many duplication-loss steps there are in a scenario
produced by such an algorithm.

We now analyse the number of steps of the scenario produced by Algorithm \ref{alg:bucket}.

\begin{prop}
The number of duplication-loss steps of a scenario produced by Algorithm \ref{alg:bucket} on a permutation of size $n$ is at most $\Theta(\frac{n}{K} \log K + \frac{n^2}{K^2})$ asymptotically. \label{prop:worst-case-upper-bound} 
\end{prop}
\begin{pf}
Suppose we are at iteration $i$ of the first loop. We have to move the $\lfloor K/2 \rfloor$ elements of $L^i$ to their positions (from $n - i \lfloor K/2 \rfloor +1$ to $n-(i-1) \lfloor K/2 \rfloor$) by duplication-loss steps of width at most $K$. The worst situation is when the elements of $L^i$ are at the begining of $\pi$. But in this case, we can move to the right the elements of $L^i$ by $\lceil K/2 \rceil$ positions at every duplication-loss step, until they reach their position. The total number of duplication-loss steps in this first process is then at most 
$$\sum_{i=1}^{\lceil \frac{n-K}{\lfloor K/2 \rfloor} \rceil} \Big\lceil
\frac{n- i \lfloor K/2 \rfloor }{\lceil K/2 \rceil} \Big\rceil  =
\Theta(\frac{n^2}{K^2}) \textrm{.} $$ 


Consider now the second loop of Algorithm \ref{alg:bucket}. In each window of size $\lfloor K/2 \rfloor$, it performs at most $\lceil \log\lfloor K/2\rfloor \rceil$ duplication-loss steps (line 7) and in the leftmost window (line 9), at most $\lceil \log K \rceil$ by the result of \cite{CCMR06}. Consequently the number of duplication-loss steps produced by the second loop is
$$\Big\lceil \frac{n-K}{\lfloor K/2 \rfloor} \Big\rceil  \lceil \log \lfloor K/2 \rfloor \rceil + \lceil \log K \rceil = \Theta(\frac{n}{K} \log K) \textrm{.}
$$
We finally get that the total number of duplication-loss steps in a scenario produced by Algorithm \ref{alg:bucket} is at most $\Theta(\frac{n}{K} \log K + \frac{n^2}{K^2})$ asymptotically in the worst case.
\end{pf}

It is easily noticed that this worst case corresponds to the \emph{reversed identity} permutation $n (n-1) \ldots 21$. This corresponds to our intuition of a worst case situation in this context.
We can also notice that $\Theta(\frac{n}{K} \log K + \frac{n^2}{K^2}) = \Theta(\frac{n^2}{K^2})$ for ``small'' values of $K$, namely as long as $K = o(\frac{n}{\log n})$. If on the contrary $\frac{n}{\log n} = o(K)$ then $\Theta(\frac{n}{K} \log K + \frac{n^2}{K^2}) = \Theta(\frac{n}{K} \log K)$. When $K = \Theta(\frac{n}{\log n})$, the two terms are of the same order.

We can also compute the average number of duplication-loss steps of a scenario produced by Algorithm \ref{alg:bucket}.

\begin{prop}
The number of duplication-loss steps of a scenario produced by Algorithm \ref{alg:bucket} on a permutation of size $n$ is on average $\Theta( \frac{n}{K} \log K + \frac{n^2}{K^2})$ asymptotically. \label{prop:average-upper-bound} 
\end{prop}
\begin{pf}
First, we introduce a few notations. Consider $\sigma$ a permutation of size $n$, and decompose it from right to left into $p=  \Big\lceil \frac{n-K}{\lfloor K/2 \rfloor} \Big\rceil +1$ windows of width $\lfloor K/2 \rfloor$, except the leftmost one, whose width is $n-\Big\lceil \frac{n-K}{\lfloor K/2 \rfloor} \Big\rceil \lfloor K/2 \rfloor \leq K$. We denote $\sigma = \sigma^1 \sigma^2 \ldots \sigma^p$ this decomposition.

Now, let us denote $c(\sigma)$ the number of duplication-loss steps produced in the first loop of Algorithm \ref{alg:bucket} on $\sigma$. And in particular, we denote $c_p(\sigma)$ the number of such steps produced by the first iteration of this loop, \emph{i.e.}, the number of steps to move the elements of $L^1$ at the end of the permutation. For computing the average number of such steps, we consider $u_n = \sum_{\sigma \in S_n} c(\sigma)$. It is simple to conceive that
\begin{eqnarray*}
u_n & = & \sum_{\sigma \in S_n} c_p(\sigma) + c(\sigma^1\ldots\sigma^{p-1}) \\
& = & \sum_{\sigma \in S_n} c_p(\sigma) + n(n-1)\ldots(n-\lfloor K/2 \rfloor+1)\sum_{\sigma \in S_{n-\lfloor K/2 \rfloor}} c(\sigma) \\
& = & \sum_{\sigma \in S_n} c_p(\sigma) + \frac{n!}{(n-\lfloor K/2 \rfloor)!} u_{n-\lfloor K/2 \rfloor} \textrm{.}
\end{eqnarray*}
Let us focus on $\sum_{\sigma \in S_n} c_p(\sigma)$. Figure \ref{fig:min} should convince the reader that
$$ \frac{n+1 -\lfloor K/2 \rfloor -\min(\sigma^p)}{K} \leq c_p(\sigma) \leq \frac{n+1-\lfloor K/2 \rfloor -\min(\sigma^p)}{\lfloor K/2 \rfloor} \textrm{.}
$$
\begin{figure}[ht]
\begin{center}
\psset{unit=0.3cm}
\begin{pspicture}(0,-3)(30,5)
\rput(1,0){$\sigma =$}
\psline(3,-0.5)(18,-0.5)
\psline(3,1)(18,1)
\psline(3,-0.5)(3,1)
\psline(18,-0.5)(18,1)
\psline(15,-0.5)(15,1)
\psline[linecolor=red]{->}(15.5,1.5)(9,1.5)
\psline[linecolor=red]{->}(16.5,2)(7,2)
\psline[linecolor=red]{->}(17.5,2.5)(13,2.5)
\rput(12,3){\color{red}$vp$-vectors}
\pscircle*(7,0.2){0.3}
\rput(0,-1.5){positions:}
\rput(3.2,-1.5){$1$}
\rput(4.,-1.5){$2$}
\rput(5.5,-1.5){$\ldots$}
\rput(8.5,-1.5){$\min(\sigma^p)$}
\rput(12,-1.5){$\ldots$}
\rput(17.5,-1.5){$n$}
\psframe[fillstyle=vlines,hatchcolor=gray,linewidth=0.08,hatchangle=-45,hatchsep=0.2](15,-0.5)(18,1)
\rput(25,-2.5){$\lfloor K/2 \rfloor$ rightmost positions}
\psline{->}(22.5,-1.5)(16,0.5)
\pcline{<->}(7,4)(18,4)
\rput(12,4.5){$n+1-\min(\sigma^p)$ elements}
\psline[linestyle=dotted](7,4)(7,0.5)
\psline[linestyle=dotted](18,4)(18,0.5)
\end{pspicture}
\caption{Bounding $c_p(\sigma)$ \label{fig:min}}
\end{center}
\end{figure}
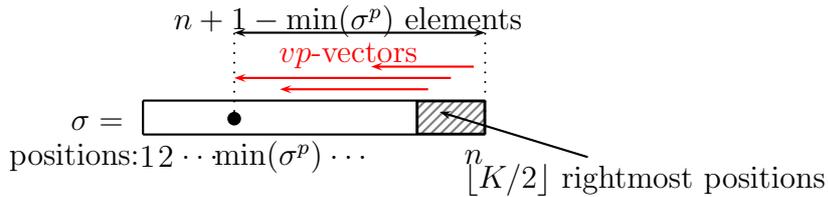

Now, we can notice that the number of permutations $\sigma$ 
of size $n$ such that $\min(\sigma^p) = i$ is 
${ n-i \choose \lfloor K/2 \rfloor -1} \left(n - \lfloor K/2 \rfloor\right)!
\lfloor K/2 \rfloor !$. This yields
\begin{eqnarray*}
&   & \sum_{\sigma \in S_n} n+1-\lfloor K/2 \rfloor-\min(\sigma^p) \\
& = & \sum_{i=1}^{n - \lfloor K/2 \rfloor +1} (n+1-\lfloor K/2 \rfloor-i) { n-i \choose \lfloor K/2 \rfloor -1} (n - \lfloor K/2 \rfloor)! \lfloor K/2 \rfloor ! \\
& = & (n - \lfloor K/2 \rfloor)! \lfloor K/2 \rfloor ! \sum_{i=\lfloor K/2 \rfloor -1}^{n - 1} (i+1-\lfloor K/2 \rfloor) { i \choose \lfloor K/2 \rfloor -1} \\
& = & (n - \lfloor K/2 \rfloor)! \lfloor K/2 \rfloor ! \lfloor K/2 \rfloor \sum_{i=\lfloor K/2 \rfloor}^{n - 1}  { i \choose \lfloor K/2 \rfloor} \\
& = & (n - \lfloor K/2 \rfloor)! \lfloor K/2 \rfloor ! \lfloor K/2 \rfloor { n \choose \lfloor K/2 \rfloor +1} \textrm{.}
\end{eqnarray*}
Consequently,
\begin{eqnarray*}
\sum_{\sigma \in S_n} c_p(\sigma) & \leq & (n - \lfloor K/2 \rfloor)! \lfloor K/2 \rfloor ! { n \choose \lfloor K/2 \rfloor +1 } \\
\sum_{\sigma \in S_n} c_p(\sigma) & \geq &  \frac{\lfloor K/2 \rfloor}{K} (n - \lfloor K/2 \rfloor)! \lfloor K/2 \rfloor ! { n \choose \lfloor K/2 \rfloor +1}\\
& \geq &  \frac{1}{3} (n - \lfloor K/2 \rfloor)! \lfloor K/2 \rfloor ! { n \choose \lfloor K/2 \rfloor +1}\textrm{,}
\end{eqnarray*}
giving after a few computations 
$$\frac{1}{3} \frac{n-\lfloor K/2 \rfloor}{\lfloor K/2 \rfloor +1} + \frac{u_{n - \lfloor K/2 \rfloor}}{(n - \lfloor K/2 \rfloor)!} \leq \frac{u_n}{n!} \leq \frac{n-\lfloor K/2 \rfloor}{\lfloor K/2 \rfloor +1} + \frac{u_{n - \lfloor K/2 \rfloor}}{(n - \lfloor K/2 \rfloor)!} \textrm{.}
$$
Therefore, we consider two sequences $(v_n)$ and $(w_n)$ satisfying the relations $v_n = \frac{1}{3} \frac{n-\lfloor K/2 \rfloor}{\lfloor K/2 \rfloor +1} + v_{n - \lfloor K/2 \rfloor}$ and $w_n = \frac{n-\lfloor K/2 \rfloor}{\lfloor K/2 \rfloor +1} + w_{n - \lfloor K/2 \rfloor}$ respectively if $n > K$, and $v_n = w_n = \frac{u_n}{n!}$ for any $n \leq K$. Then we have $v_n \leq \frac{u_n}{n!} \leq w_n \forall n \in \mathbb{N}$. We can solve the recurrence equations for $v_n$ and $w_n$; and if we write $n = \lceil \frac{n-K}{\lfloor K/2 \rfloor} \rceil \lfloor K/2 \rfloor +r$ (then $\lfloor K/2 \rfloor \leq r \leq K$), we get:
\begin{eqnarray*}
v_n & = & \frac{1}{3} \sum_{i=1}^{\lceil \frac{n-K}{\lfloor K/2 \rfloor} \rceil} \frac{n - i \lfloor K/2 \rfloor}{\lfloor K/2 \rfloor +1} + v_r \\
& = & \frac{1}{3 (\lfloor K/2 \rfloor +1)} \Big\lceil \frac{n-K}{\lfloor K/2 \rfloor} \Big \rceil \Big( n - \lfloor K/2 \rfloor \frac{\lceil \frac{n-K}{\lfloor K/2 \rfloor} \rceil +1}{2}\Big) + v_r \\
& = & \Theta(\frac{n^2}{K^2})
\end{eqnarray*}
and 
$$w_n  =   \sum_{i=1}^{\lceil \frac{n-K}{\lfloor K/2 \rfloor} \rceil} \frac{n
- i \lfloor K/2 \rfloor}{\lfloor K/2 \rfloor +1} + w_r  =
\Theta(\frac{n^2}{K^2})$$

Consequently, the average number of duplication-loss steps produced by the first loop of Algorithm \ref{alg:bucket} on permutations of size $n$ is $\frac{u_n}{n!} = \Theta(\frac{n^2}{K^2})$.

What is left to compute is the average number of duplication-loss steps produced by the second loop of Algorithm \ref{alg:bucket} on permutations of size $n$. This number is given by
\begin{eqnarray*}
& & \frac{1}{n!} \sum_{\sigma \in S_n} \sum_{i=1}^p \lceil \log(desc(\sigma^i)+1) \rceil \\
& = & \frac{1}{n!} \Big( \sum_{i=2}^p \sum_{\sigma \in S_n} \lceil \log(desc(\sigma^i)+1) \rceil + \sum_{\sigma \in S_n} \lceil \log(desc(\sigma^1)+1) \rceil \Big) \\
& = & \frac{1}{n!} \Big( \sum_{i=2}^p (n- \lfloor K/2 \rfloor) ! { n \choose \lfloor K/2 \rfloor } \sum_{\sigma \in S_{\lfloor K/2 \rfloor}} \lceil \log(desc(\sigma)+1) \rceil  \\
& & + (n- |\sigma^1|) ! { n \choose |\sigma^1|} \sum_{\sigma \in S_{|\sigma^1|}} \lceil \log(desc(\sigma)+1) \rceil \Big) \\
& = & \frac{1}{\lfloor K/2 \rfloor !} (p-1) \sum_{\sigma \in S_{\lfloor K/2 \rfloor}} \lceil \log(desc(\sigma)+1) \rceil \\
& & + \frac{1}{|\sigma^1|!} \sum_{\sigma \in S_{|\sigma^1|}} \lceil \log(desc(\sigma)+1) \rceil \textrm{.}
\end{eqnarray*}

Since $p= \Big\lceil \frac{n-K}{\lfloor K/2 \rfloor} \Big \rceil +1$, we deduce that the average number of duplication-loss steps produced by the second loop of Algorithm \ref{alg:bucket} on permutations of size $n$ is $\Theta(\frac{1}{\lfloor K/2 \rfloor !}  \Big( \Big\lceil \frac{n-K}{\lfloor K/2 \rfloor} \Big \rceil +1 \Big) \sum_{\sigma \in S_{\lfloor K/2 \rfloor}} \lceil \log(desc(\sigma)+1) \rceil)$. Hence we focus on the computation of $\frac{1}{k!}\sum_{\sigma \in S_{k}} \lceil \log(desc(\sigma)+1) \rceil)$ for $k= \lfloor K/2 \rfloor$.
By concavity of the $\log$ function, since $\frac{1}{k!} \sum_{\sigma \in S_k} desc(\sigma)+1 = \frac{k+1}{2}$, we get that $$\frac{1}{k!}\sum_{\sigma \in S_{k}} \lceil \log(desc(\sigma)+1) \rceil \geq \frac{1}{k!}\sum_{\sigma \in S_{k}} \log(desc(\sigma)+1) \geq \log(\frac{k+1}{2}) \textrm{.}$$

Moreover, it is clear that
$$
\frac{1}{k!}\sum_{\sigma \in S_{k}} \lceil \log(desc(\sigma)+1) \rceil \leq \lceil \log(k) \rceil \textrm{,}
$$
so that we deduce that 
\begin{equation} \label{eq:log-average}
\frac{1}{k!}\sum_{\sigma \in S_{k}} \lceil \log(desc(\sigma)+1) \rceil = \Theta ( \log(k) ) \textrm{.}
\end{equation}

Consequently, the average number of duplication-loss steps produced by the second loop of Algorithm \ref{alg:bucket} on permutations of size $n$ is $\Theta ( \lceil \frac{n-K}{\lfloor K/2 \rfloor} \rceil \log(\lfloor K/2 \rfloor) ) = \Theta (\frac{n}{K} \log K) $.

Finally, we end the proof concluding that the total number of duplication-loss steps in a scenario produced by Algorithm \ref{alg:bucket} on a permutation of size $n$ is $\Theta (\frac{n}{K} \log K + \frac{n^2}{K^2}) $ on average.
\end{pf}

\subsection{Lower bound}

It is possible to provide very simple lower bounds on the number of
duplication-loss steps necessary to obtain a permutation. These lower bounds are given and proved in Propositions \ref{prop:worst-case-lower-bound} and \ref{prop:average-lower-bound} below. They are tight in most cases, however not in any case. Indeed the upper and lower bounds coincide up to a constant factor whenever $K$ is a constant, or when $K = K(n)$, except when $\frac{n}{ \log n} \ll K(n) \ll n$.

\begin{prop}
In the worst case, $\Omega(\log n + \frac{n^2}{K^2})$ duplication-loss steps of width $K$ are necessary to obtain a permutation of $S_n$ from $123\ldots n$.
\label{prop:worst-case-lower-bound}
\end{prop}
\begin{pf}
Let us consider first the number of inversions in a permutation that can create a duplication-loss step $s$ of width $K$. It is easily seen that these new inversions can only involve two elements of $s$. Call $i$ the number of elements of $s$ that are kept in the first copy. Then the maximum number of inversions that can be created by $s$ is $i(K-i) \leq \frac{K^2}{4}$. Now, a permutation $\sigma \in S_n$ has up to $\frac{n(n-1)}{2}$ inversions, so that at least $\frac{2n(n-1)}{K^2}$ duplication-loss steps are necessary to transform $123\ldots n$ into $\sigma$.

To get the other term of the lower bound, we just refer to the result of \cite{CCMR06} recalled at the beginning of this section, namely that $\log n$ steps are necessary in the worst case in the whole genome duplication - random loss model, in which duplication-loss operation are less restricted.

Finally, we get a lower bound of $\Omega( \log n + \frac{n^2}{K^2})$ necessary duplication-loss steps to obtain a permutation of $S_n$ from $123 \ldots n$ in the worst case.
\end{pf}

\begin{prop}
On average, $\Omega( \log n + \frac{n^2}{K^2})$ duplication-loss steps of width $K$ are necessary to obtain a permutation of $S_n$ from $123\ldots n$.
\label{prop:average-lower-bound}
\end{prop}
\begin{pf}
As before, a duplication-loss step can create at most $\frac{K^2}{4}$ inversions in a permutation. But the average number of inversions in a permutation of $S_n$ is $\frac{n(n-1)}{4}$, so that on average at least $\frac{n(n-1)}{K^2}$ duplication-loss steps are necessary to transform $123\ldots n$ into $\sigma \in S_n$.

Again, \cite{CCMR06} provides use with the $\Omega(\log n)$ lower bound,
referring to the whole genome duplication - random loss model which is more general than ours, so that this bound applies in our context.

We conclude that a lower bound on the average number of duplication-loss steps necessary to obtain a permutation of $S_n$ from $123 \ldots n$ is $\Omega( \log n + \frac{n^2}{K^2})$.
\end{pf}

\section{Conclusion}

We discuss the results of Section \ref{section:algo} on the average (or worst case) number of steps of width $K$ to obtain a permutation of size $n$. It appears that we could not provide lower bounds that coincide with the upper bounds given by our algorithm, but we claim that they are tight in many cases. Indeed, whenever $K =
o(\frac{n}{\log n})$, we get that $\frac{n}{K} \log K = o(\frac{n^2}{K^2})$, and
consequently the upper bound can be rewritten as $\Theta(\frac{n}{K} \log K  +
\frac{n^2}{K^2}) = \Theta(\frac{n^2}{K^2})$, which coincide up to a constant
factor with the lower bound $\Omega(\log n  + \frac{n^2}{K^2}) =
\Omega(\frac{n^2}{K^2})$. For the case $K = \Theta(\frac{n}{\log n})$, the same
argument holds, but the constant factor between the lower and the upper bound
might be much greater. Finally, if $K = \Theta(n)$, then $\Theta(\frac{n}{K} \log
K  + \frac{n^2}{K^2}) = \Theta(\log n)$ and $\Omega(\log n  + \frac{n^2}{K^2}) =
\Omega (\log n)$, so that upper and lower bounds coincide again.

On the contrary, when $ \frac{n}{\log n} \ll K \ll n$, the upper and lower bounds
provided do not coincide. We leave as an open question the problem of finding an
algorithm that computes a duplication-loss scenario whose number of steps is
optimal (on average and in the worst case) up to a constant factor, when the
width $K$ of the duplicated windows satisfies $ \frac{n}{\log n} \ll K \ll n$.

Several other questions are still open. First of all neither of our algorithms is
optimal for a specific permutation and our results are only optimal asymptotically
in average and/or in the worst case. It could be interesting to provide
algorithms that produce optimal duplication-loss scenarios on any permutation $\sigma$, for $K = K(n)$ in order to provide some continuity between the
bubble sort (corresponding to $K=2$) and the radix sort (corresponding to $K(n)=n$).
\nocite{*}



\end{document}